

\input amstex
\expandafter\ifx\csname mathdefs.tex\endcsname\relax
  \expandafter\gdef\csname mathdefs.tex\endcsname{}
\else \message{Hey!  Apparently you were trying to
  \string twice.   This does not make sense.} 
\errmessage{Please edit your file (probably \jobname.tex) and remove
any duplicate ``\string\input'' lines} \fi




\catcode`\X=12\catcode`\@=11

\def\n@wcount{\alloc@0\count\countdef\insc@unt}
\def\n@wwrite{\alloc@7\write\chardef\sixt@@n}
\def\n@wread{\alloc@6\read\chardef\sixt@@n}
\def\r@s@t{\relax}\def\v@idline{\par}\def\@mputate#1/{#1}
\def\l@c@l#1X{\firstpart.#1}\def\gl@b@l#1X{#1}\def\t@d@l#1X{{}}

\def\crossrefs#1{\ifx\all#1\let\tr@ce=\all\else\def\tr@ce{#1,}\fi
   \n@wwrite\cit@tionsout\openout\cit@tionsout=\jobname.cit 
   \write\cit@tionsout{\tr@ce}\expandafter\setfl@gs\tr@ce,}
\def\setfl@gs#1,{\def\@{#1}\ifx\@\empty\let\next=\relax
   \else\let\next=\setfl@gs\expandafter\xdef
   \csname#1tr@cetrue\endcsname{}\fi\next}
\def\m@ketag#1#2{\expandafter\n@wcount\csname#2tagno\endcsname
     \csname#2tagno\endcsname=0\let\tail=\all\xdef\all{\tail#2,}
   \ifx#1\l@c@l\let\tail=\r@s@t\xdef\r@s@t{\csname#2tagno\endcsname=0\tail}\fi
   \expandafter\gdef\csname#2cite\endcsname##1{\expandafter
     \ifx\csname#2tag##1\endcsname\relax?\else\csname#2tag##1\endcsname\fi
     \expandafter\ifx\csname#2tr@cetrue\endcsname\relax\else
     \write\cit@tionsout{#2tag ##1 cited on page \folio.}\fi}
   \expandafter\gdef\csname#2page\endcsname##1{\expandafter
     \ifx\csname#2page##1\endcsname\relax?\else\csname#2page##1\endcsname\fi
     \expandafter\ifx\csname#2tr@cetrue\endcsname\relax\else
     \write\cit@tionsout{#2tag ##1 cited on page \folio.}\fi}
   \expandafter\gdef\csname#2tag\endcsname##1{\expandafter
      \ifx\csname#2check##1\endcsname\relax
      \expandafter\xdef\csname#2check##1\endcsname{}%
      \else\immediate\write16{Warning: #2tag ##1 used more than once.}\fi
      \multit@g{#1}{#2}##1/X%
      \write\t@gsout{#2tag ##1 assigned number \csname#2tag##1\endcsname\space
      on page \number\count0.}%
   \csname#2tag##1\endcsname}}
\def\multit@g#1#2#3/#4X{\def\t@mp{#4}\ifx\t@mp\empty%
      \global\advance\csname#2tagno\endcsname by 1 
      \expandafter\xdef\csname#2tag#3\endcsname
      {#1\number\csname#2tagno\endcsnameX}%
   \else\expandafter\ifx\csname#2last#3\endcsname\relax
      \expandafter\n@wcount\csname#2last#3\endcsname
      \global\advance\csname#2tagno\endcsname by 1 
      \expandafter\xdef\csname#2tag#3\endcsname
      {#1\number\csname#2tagno\endcsnameX}
      \write\t@gsout{#2tag #3 assigned number \csname#2tag#3\endcsname\space
      on page \number\count0.}\fi
   \global\advance\csname#2last#3\endcsname by 1
   \def\t@mp{\expandafter\xdef\csname#2tag#3/}%
   \expandafter\t@mp\@mputate#4\endcsname
   {\csname#2tag#3\endcsname\lastpart{\csname#2last#3\endcsname}}\fi}
\def\t@gs#1{\def\all{}\m@ketag#1e\m@ketag#1s\m@ketag\t@d@l p
   \m@ketag\gl@b@l r \n@wread\t@gsin
   \openin\t@gsin=\jobname.tgs \re@der \closein\t@gsin
   \n@wwrite\t@gsout\openout\t@gsout=\jobname.tgs }
\outer\def\localtags{\t@gs\l@c@l}
\outer\def\globaltags{\t@gs\gl@b@l}
\outer\def\newlocaltag#1{\m@ketag\l@c@l{#1}}
\outer\def\newglobaltag#1{\m@ketag\gl@b@l{#1}}

\newif\ifpr@ 
\def\m@kecs #1tag #2 assigned number #3 on page #4.%
   {\expandafter\gdef\csname#1tag#2\endcsname{#3}
   \expandafter\gdef\csname#1page#2\endcsname{#4}
   \ifpr@\expandafter\xdef\csname#1check#2\endcsname{}\fi}
\def\re@der{\ifeof\t@gsin\let\next=\relax\else
   \read\t@gsin to\t@gline\ifx\t@gline\v@idline\else
   \expandafter\m@kecs \t@gline\fi\let \next=\re@der\fi\next}
\def\pretags#1{\pr@true\pret@gs#1,,}
\def\pret@gs#1,{\def\@{#1}\ifx\@\empty\let\n@xtfile=\relax
   \else\let\n@xtfile=\pret@gs \openin\t@gsin=#1.tgs \message{#1} \re@der 
   \closein\t@gsin\fi \n@xtfile}

\newcount\sectno\sectno=0\newcount\subsectno\subsectno=0
\newif\ifultr@local \def\ultralocal{\ultr@localtrue}
\def\firstpart{\number\sectno}
\def\lastpart#1{\ifcase#1 \or a\or b\or c\or d\or e\or f\or g\or h\or 
   i\or k\or l\or m\or n\or o\or p\or q\or r\or s\or t\or u\or v\or w\or 
   x\or y\or z \fi}

\def\resetall{\global\advance\sectno by 1\subsectno=0
   \gdef\firstpart{\number\sectno}\r@s@t}
\def\resetsub{\global\advance\subsectno by 1
   \gdef\firstpart{\number\sectno.\number\subsectno}\r@s@t}
\def\newsection#1\par{\resetall\vskip0pt plus.3\vsize\penalty-250
   \vskip0pt plus-.3\vsize\bigskip\bigskip
   \message{#1}\leftline{\bf#1}\nobreak\bigskip}
\def\subsection#1\par{\ifultr@local\resetsub\fi
   \vskip0pt plus.2\vsize\penalty-250\vskip0pt plus-.2\vsize
   \bigskip\smallskip\message{#1}\leftline{\bf#1}\nobreak\medskip}

\def\t@gsoff#1,{\def\@{#1}\ifx\@\empty\let\next=\relax\else\let\next=\t@gsoff
   \def\@@{p}\ifx\@\@@\else
   \expandafter\gdef\csname#1cite\endcsname##1{\zeigen{##1}}
   \expandafter\gdef\csname#1page\endcsname##1{?}
   \expandafter\gdef\csname#1tag\endcsname##1{\zeigen{##1}}\fi\fi\next}
\def\verbatimtags{\ifx\all\relax\else\expandafter\t@gsoff\all,\fi}
\def\zeigen#1{\hbox{$\langle$}#1\hbox{$\rangle$}}

\def\(#1){\edef\dot@g{\ifmmode\ifinner(\hbox{\noexpand\etag{#1}})
   \else\noexpand\eqno(\hbox{\noexpand\etag{#1}})\fi
   \else(\noexpand\ecite{#1})\fi}\dot@g}

\newif\ifbr@ck
\def\eat#1{}
\def\[#1]{\br@cktrue[\br@cket#1'X]}
\def\br@cket#1'#2X{\def\temp{#2}\ifx\temp\empty\let\next\eat
   \else\let\next\br@cket\fi
   \ifbr@ck\br@ckfalse\br@ck@t#1,X\else\br@cktrue#1\fi\next#2X}
\def\br@ck@t#1,#2X{\def\temp{#2}\ifx\temp\empty\let\neext\eat
   \else\let\neext\br@ck@t\def\temp{,}\fi
   \def\teemp{#1}\ifx\teemp\empty\else\rcite{#1}\fi\temp\neext#2X}
\def\resetbr@cket{\gdef\[##1]{[\rtag{##1}]}}
\def\references{\resetbr@cket\newsection References\par}

\newtoks\symb@ls\newtoks\s@mb@ls\newtoks\p@gelist\n@wcount\ftn@mber
    \ftn@mber=1\newif\ifftn@mbers\ftn@mbersfalse\newif\ifbyp@ge\byp@gefalse
\def\defm@rk{\ifftn@mbers\n@mberm@rk\else\symb@lm@rk\fi}
\def\n@mberm@rk{\xdef\m@rk{{\the\ftn@mber}}%
    \global\advance\ftn@mber by 1 }
\def\rot@te#1{\let\temp=#1\global#1=\expandafter\r@t@te\the\temp,X}
\def\r@t@te#1,#2X{{#2#1}\xdef\m@rk{{#1}}}
\def\b@@st#1{{$^{#1}$}}\def\str@p#1{#1}
\def\symb@lm@rk{\ifbyp@ge\rot@te\p@gelist\ifnum\expandafter\str@p\m@rk=1 
    \s@mb@ls=\symb@ls\fi\write\f@nsout{\number\count0}\fi \rot@te\s@mb@ls}
\def\byp@ge{\byp@getrue\n@wwrite\f@nsin\openin\f@nsin=\jobname.fns 
    \n@wcount\currentp@ge\currentp@ge=0\p@gelist={0}
    \re@dfns\closein\f@nsin\rot@te\p@gelist
    \n@wread\f@nsout\openout\f@nsout=\jobname.fns }
\def\m@kelist#1X#2{{#1,#2}}
\def\re@dfns{\ifeof\f@nsin\let\next=\relax\else\read\f@nsin to \f@nline
    \ifx\f@nline\v@idline\else\let\t@mplist=\p@gelist
    \ifnum\currentp@ge=\f@nline
    \global\p@gelist=\expandafter\m@kelist\the\t@mplistX0
    \else\currentp@ge=\f@nline
    \global\p@gelist=\expandafter\m@kelist\the\t@mplistX1\fi\fi
    \let\next=\re@dfns\fi\next}
\def\symbols#1{\symb@ls={#1}\s@mb@ls=\symb@ls} 
\def\bigsymbol{\textstyle}
\symbols{\bigsymbol\ast,\dagger,\ddagger,\sharp,\flat,\natural,\star}
\def\ftnumbers{\ftn@mberstrue} \def\ftsymbols{\ftn@mbersfalse}
\def\paginal{\byp@ge} \def\resetftnumbers{\ftn@mber=1}
\def\ftnote#1{\defm@rk\expandafter\expandafter\expandafter\footnote
    \expandafter\b@@st\m@rk{#1}}

\long\def\jump#1\endjump{}
\def\ssum{\mathop{\lower .1em\hbox{$\textstyle\Sigma$}}\nolimits}

\def\qed{\nobreak\kern 1em \vrule height .5em width .5em depth 0em}
\def\newneq{\hbox{\rlap{\hbox to 1\wd9{\hss$=$\hss}}\raise .1em 
   \hbox to 1\wd9{\hss$\scriptscriptstyle/$\hss}}}
\def\subsetne{\setbox9 = \hbox{$\subset$}\mathrel{\hbox{\rlap
   {\lower .4em \newneq}\raise .13em \hbox{$\subset$}}}}
\def\supsetne{\setbox9 = \hbox{$\subset$}\mathrel{\hbox{\rlap
   {\lower .4em \newneq}\raise .13em \hbox{$\supset$}}}}

\def\vbar{\mathchoice{\vrule height6.3ptdepth-.5ptwidth.8pt\kern-.8pt}
   {\vrule height6.3ptdepth-.5ptwidth.8pt\kern-.8pt}
   {\vrule height4.1ptdepth-.35ptwidth.6pt\kern-.6pt}
   {\vrule height3.1ptdepth-.25ptwidth.5pt\kern-.5pt}}
\def\f@dge{\mathchoice{}{}{\mkern.5mu}{\mkern.8mu}}
\def\b@c#1#2{{\rm \mkern#2mu\vbar\mkern-#2mu#1}}
\def\b@b#1{{\rm I\mkern-3.5mu #1}}
\def\b@a#1#2{{\rm #1\mkern-#2mu\f@dge #1}}
\def\bb#1{{\count4=`#1 \advance\count4by-64 \ifcase\count4\or\b@a A{11.5}\or
   \b@b B\or\b@c C{5}\or\b@b D\or\b@b E\or\b@b F \or\b@c G{5}\or\b@b H\or
   \b@b I\or\b@c J{3}\or\b@b K\or\b@b L \or\b@b M\or\b@b N\or\b@c O{5} \or
   \b@b P\or\b@c Q{5}\or\b@b R\or\b@a S{8}\or\b@a T{10.5}\or\b@c U{5}\or
   \b@a V{12}\or\b@a W{16.5}\or\b@a X{11}\or\b@a Y{11.7}\or\b@a Z{7.5}\fi}}

\catcode`\X=11 \catcode`\@=12

\expandafter\ifx\csname citeadd.tex\endcsname\relax
\expandafter\gdef\csname citeadd.tex\endcsname{}
\else \message{Hey!  Apparently you were trying to
\string twice.   This does not make sense.} 
\errmessage{Please edit your file (probably \jobname.tex) and remove
any duplicate ``\string\input'' lines} \fi

\def\sciteu{\sciteerror{undefined}}

\def\sciteerror#1#2{{\mathortextbf{\scite{#2}}}\complainaboutcitation{#1}{#2}}
\def\mathortextbf#1{\hbox{\bf #1}}
\def\complainaboutcitation#1#2{%
\vadjust{\line{\llap{---$\!\!>$ }\qquad scite$\{$#2$\}$ #1\hfil}}}

\sectno=-1   
\localtags
\NoBlackBoxes
\ifx\shlhetal\undefinedcontrolsequence\let\shlhetal\relax\fi
\define\mr{\medskip\roster}
\define\sn{\smallskip\noindent}
\define\mn{\medskip\noindent}
\define\bn{\bigskip\noindent}
\define\ub{\underbar}
\define\wilog{\text{without loss of generality}}
\define\ermn{\endroster\medskip\noindent}

\define\dbcu{\dsize\bigcup}
\define\nl{\newline}
\documentstyle {amsppt}
\topmatter
\title{Polynomial Time Logic: Inability to express \\
 Sh634} \endtitle
\author {Saharon Shelah \thanks {\null\newline 
Partially supported by the United States-Israel Binational Science
Foundation \null\newline
I would like to thank Alice Leonhardt for the beautiful typing. \null\newline
 Done - Fall 1996 \null\newline 
 Latest Revision - 98/July/23} \endthanks} \endauthor 
\affil{Institute of Mathematics\\
 The Hebrew University\\
 Jerusalem, Israel
 \medskip
 Rutgers University\\
 Mathematics Department\\
 New Brunswick, NJ  USA} \endaffil
\endtopmatter
\document  

\expandafter\ifx\csname alice2jlem.tex\endcsname\relax
  \expandafter\gdef\csname alice2jlem.tex\endcsname{}
\else \message{Hey!  Apparently you were trying to
\string  twice.   This does not make sense.}
\errmessage{Please edit your file (probably \jobname.tex) and remove
any duplicate ``\string\input'' lines} \fi

\expandafter\ifx\csname bib4plain.tex\endcsname\relax
  \expandafter\gdef\csname bib4plain.tex\endcsname{}
\else \message{Hey!  Apparently you were trying to \string twice.   This does not make sense.}
\errmessage{Please edit your file (probably \jobname.tex) and remove
any duplicate ``\string\input'' lines} \fi

\def\renewcommand{\newcommand}	       
\edef\cite{\the\catcode`@}%
\catcode`@ = 11
\let\@oldatcatcode = \cite
\chardef\@letter = 11
\chardef\@other = 12
%
%
%
%
\def\@innerdef#1#2{\edef#1{\expandafter\noexpand\csname #2\endcsname}}%
%
%
\@innerdef\@innernewcount{newcount}%
\@innerdef\@innernewdimen{newdimen}%
\@innerdef\@innernewif{newif}%
\@innerdef\@innernewwrite{newwrite}%
%
%
%
\def\@gobble#1{}%
%
%
%
\ifx\inputlineno\@undefined
   \let\@linenumber = \empty 
\else
   \def\@linenumber{\the\inputlineno:\space}%
\fi
%
%
%
\def\@futurenonspacelet#1{\def\cs{#1}%
   \afterassignment\@stepone\let\@nexttoken=
}%
\begingroup 
\def\\{\global\let\@stoken= }%
\\ 
\endgroup
\def\@stepone{\expandafter\futurelet\cs\@steptwo}%
\def\@steptwo{\expandafter\ifx\cs\@stoken\let\@@next=\@stepthree
   \else\let\@@next=\@nexttoken\fi \@@next}%
\def\@stepthree{\afterassignment\@stepone\let\@@next= }%
%
%
%
\def\@getoptionalarg#1{%
   \let\@optionaltemp = #1%
   \let\@optionalnext = \relax
   \@futurenonspacelet\@optionalnext\@bracketcheck
}%
%
%
\def\@bracketcheck{%
   \ifx [\@optionalnext
      \expandafter\@@getoptionalarg
   \else
      \let\@optionalarg = \empty
      \expandafter\@optionaltemp
   \fi
}%
\def\@@getoptionalarg[#1]{%
   \def\@optionalarg{#1}%
   \@optionaltemp
}%
%
%
%
\def\@nnil{\@nil}%
\def\@fornoop#1\@@#2#3{}%
\def\@for#1:=#2\do#3{%
   \edef\@fortmp{#2}%
   \ifx\@fortmp\empty \else
      \expandafter\@forloop#2,\@nil,\@nil\@@#1{#3}%
   \fi
}%
\def\@forloop#1,#2,#3\@@#4#5{\def#4{#1}\ifx #4\@nnil \else
       #5\def#4{#2}\ifx #4\@nnil \else#5\@iforloop #3\@@#4{#5}\fi\fi
}%
\def\@iforloop#1,#2\@@#3#4{\def#3{#1}\ifx #3\@nnil
       \let\@nextwhile=\@fornoop \else
      #4\relax\let\@nextwhile=\@iforloop\fi\@nextwhile#2\@@#3{#4}%
}%
%
%
%
\@innernewif\if@fileexists
\def\@testfileexistence{\@getoptionalarg\@finishtestfileexistence}%
\def\@finishtestfileexistence#1{%
   \begingroup
      \def\extension{#1}%
      \immediate\openin0 =
         \ifx\@optionalarg\empty\jobname\else\@optionalarg\fi
         \ifx\extension\empty \else .#1\fi
         \space
      \ifeof 0
         \global\@fileexistsfalse
      \else
         \global\@fileexiststrue
      \fi
      \immediate\closein0
   \endgroup
}%
%
%
%
%
\def\bibliographystyle#1{%
   \@readauxfile
   \@writeaux{\string\bibstyle{#1}}%
}%
\let\bibstyle = \@gobble
%
%
\let\bblfilebasename = \jobname
\def\bibliography#1{%
   \@readauxfile
   \@writeaux{\string\bibdata{#1}}%
   \@testfileexistence[\bblfilebasename]{bbl}%
   \if@fileexists
      \nobreak
      \@readbblfile
   \fi
}%
\let\bibdata = \@gobble
%
%
\def\nocite#1{%
   \@readauxfile
   \@writeaux{\string\citation{#1}}%
}%
\@innernewif\if@notfirstcitation
%
%
\def\cite{\@getoptionalarg\@cite}%
%
%
\def\@cite#1{%
   \let\@citenotetext = \@optionalarg
   \printcitestart
   \nocite{#1}%
   \@notfirstcitationfalse
   \@for \@citation :=#1\do
   {%
      \expandafter\@onecitation\@citation\@@
   }%
   \ifx\empty\@citenotetext\else
      \printcitenote{\@citenotetext}%
   \fi
   \printcitefinish
}%
\def\@onecitation#1\@@{%
   \if@notfirstcitation
      \printbetweencitations
   \fi
   \expandafter \ifx \csname\@citelabel{#1}\endcsname \relax
      \if@citewarning
         \message{\@linenumber Undefined citation `#1'.}%
      \fi
      \expandafter\gdef\csname\@citelabel{#1}\endcsname{%
\strut
\vadjust{\vskip-\dp\strutbox
\vbox to 0pt{\vss\parindent0cm \leftskip=\hsize 
\advance\leftskip3mm
\advance\hsize 4cm\strut\openup-4pt 
\rightskip 0cm plus 1cm minus 0.5cm ?  #1 ?\strut}}
         {\tt
            \escapechar = -1
            \nobreak\hskip0pt
            \expandafter\string\csname#1\endcsname
            \nobreak\hskip0pt
         }%
      }%
   \fi
   \csname\@citelabel{#1}\endcsname
   \@notfirstcitationtrue
}%
%
%
\def\@citelabel#1{b@#1}%
%
%
\def\@citedef#1#2{\expandafter\gdef\csname\@citelabel{#1}\endcsname{#2}}%
%
%
%
\def\@readbblfile{%
   \ifx\@itemnum\@undefined
      \@innernewcount\@itemnum
   \fi
   \begingroup
      \def\begin##1##2{%
         \setbox0 = \hbox{\biblabelcontents{##2}}%
         \biblabelwidth = \wd0
      }%
      \def\end##1{}
      %
      %
      \@itemnum = 0
      \def\bibitem{\@getoptionalarg\@bibitem}%
      \def\@bibitem{%
         \ifx\@optionalarg\empty
            \expandafter\@numberedbibitem
         \else
            \expandafter\@alphabibitem
         \fi
      }%
      \def\@alphabibitem##1{%
         \expandafter \xdef\csname\@citelabel{##1}\endcsname {\@optionalarg}%
         \ifx\biblabelprecontents\@undefined
            \let\biblabelprecontents = \relax
         \fi
         \ifx\biblabelpostcontents\@undefined
            \let\biblabelpostcontents = \hss
         \fi
         \@finishbibitem{##1}%
      }%
      \def\@numberedbibitem##1{%
         \advance\@itemnum by 1
         \expandafter \xdef\csname\@citelabel{##1}\endcsname{\number\@itemnum}%
         \ifx\biblabelprecontents\@undefined
            \let\biblabelprecontents = \hss
         \fi
         \ifx\biblabelpostcontents\@undefined
            \let\biblabelpostcontents = \relax
         \fi
         \@finishbibitem{##1}%
      }%
      \def\@finishbibitem##1{%
         \biblabelprint{\csname\@citelabel{##1}\endcsname}%
         \@writeaux{\string\@citedef{##1}{\csname\@citelabel{##1}\endcsname}}%
         \ignorespaces
      }%
      %
      %
      \let\em = \bblem
      \let\newblock = \bblnewblock
      \let\sc = \bblsc
      \frenchspacing
      \clubpenalty = 4000 \widowpenalty = 4000
      \tolerance = 10000 \hfuzz = .5pt
      \everypar = {\hangindent = \biblabelwidth
                      \advance\hangindent by \biblabelextraspace}%
      \bblrm
      \parskip = 1.5ex plus .5ex minus .5ex
      \biblabelextraspace = .5em
      \bblhook
      \input \bblfilebasename.bbl
   \endgroup
}%
%
%
\@innernewdimen\biblabelwidth
\@innernewdimen\biblabelextraspace
%
%
%
\def\biblabelprint#1{%
   \noindent
   \hbox to \biblabelwidth{%
      \biblabelprecontents
      \biblabelcontents{#1}%
      \biblabelpostcontents
   }%
   \kern\biblabelextraspace
}%
%
%
%
\def\biblabelcontents#1{{\bblrm [#1]}}%
%
%
\def\bblrm{\rm}%
%
%
\def\bblem{\it}%
%
%
\def\bblsc{\ifx\@scfont\@undefined
              \font\@scfont = cmcsc10
           \fi
           \@scfont
}%
%
%
\def\bblnewblock{\hskip .11em plus .33em minus .07em }%
%
%
\let\bblhook = \empty
%
%
%
\def\printcitestart{[}
\def\printcitefinish{]}
\def\printbetweencitations{, }
\def\printcitenote#1{, #1}
%
%
%
\let\citation = \@gobble
%
%
%
\@innernewcount\@numparams
%
%
\def\newcommand#1{%
   \def\@commandname{#1}%
   \@getoptionalarg\@continuenewcommand
}%
%
%
\def\@continuenewcommand{%
   \@numparams = \ifx\@optionalarg\empty 0\else\@optionalarg \fi \relax
   \@newcommand
}%
%
%
\def\@newcommand#1{%
   \def\@startdef{\expandafter\edef\@commandname}%
   \ifnum\@numparams=0
      \let\@paramdef = \empty
   \else
      \ifnum\@numparams>9
         \errmessage{\the\@numparams\space is too many parameters}%
      \else
         \ifnum\@numparams<0
            \errmessage{\the\@numparams\space is too few parameters}%
         \else
            \edef\@paramdef{%
               \ifcase\@numparams
                  \empty  No arguments.
               \or ####1%
               \or ####1####2%
               \or ####1####2####3%
               \or ####1####2####3####4%
               \or ####1####2####3####4####5%
               \or ####1####2####3####4####5####6%
               \or ####1####2####3####4####5####6####7%
               \or ####1####2####3####4####5####6####7####8%
               \or ####1####2####3####4####5####6####7####8####9%
               \fi
            }%
         \fi
      \fi
   \fi
   \expandafter\@startdef\@paramdef{#1}%
}%
%
%
%
%
\def\@readauxfile{%
   \if@auxfiledone \else 
      \global\@auxfiledonetrue
      \@testfileexistence{aux}%
      \if@fileexists
         \begingroup
            \endlinechar = -1
            \catcode`@ = 11
            \input \jobname.aux
         \endgroup
      \else
         \message{\@undefinedmessage}%
         \global\@citewarningfalse
      \fi
      \immediate\openout\@auxfile = \jobname.aux
   \fi
}%
%
%
\newif\if@auxfiledone
\ifx\noauxfile\@undefined \else \@auxfiledonetrue\fi
%
%
%
%
\@innernewwrite\@auxfile
\def\@writeaux#1{\ifx\noauxfile\@undefined \write\@auxfile{#1}\fi}%
%
%
%
\ifx\@undefinedmessage\@undefined
   \def\@undefinedmessage{No .aux file; I won't give you warnings about
                          undefined citations.}%
\fi
%
%
\@innernewif\if@citewarning
\ifx\noauxfile\@undefined \@citewarningtrue\fi
%
%
%
\catcode`@ = \@oldatcatcode


\def\widestnumber#1#2{}

\def\rm{\fam0 \tenrm}

\def\fakesubhead#1\endsubhead{\bigskip\noindent{\bf#1}\par}


%
%
%

%

\font\textrsfs=rsfs10
\font\scriptrsfs=rsfs7
\font\scriptscriptrsfs=rsfs5

\newfam\rsfsfam
\textfont\rsfsfam=\textrsfs
\scriptfont\rsfsfam=\scriptrsfs
\scriptscriptfont\rsfsfam=\scriptscriptrsfs

\edef\oldcatcodeofat{\the\catcode`\@}
\catcode`\@11

\def\Cal@@#1{\noaccents@ \fam \rsfsfam #1}

\catcode`\@\oldcatcodeofat

\newpage

\head {Anotated Content} \endhead  \resetall
\bn
\S1 $\quad$ The polynomial time logic presented 
\roster
\item "{{}}"  [We present this logic from a paper of Blass, Gurevich,
Shelah \cite{BGSh:533}; to compute what you can compute from a model in
polynomial time without arbitrary choices (like ordering the model).]
\endroster
\bn
\S2 $\quad$ The general elimination of quantifiers and proof it's 
non-expressive 
\roster
\item "{{}}"  [We define a criterion for showing the logic cannot say too
complicated things on some model using a family of partial automorphism 
(rather than real automorphisms) and prove that it works.  This is a relative
of the Ehrenfeucht-Fraisse games, and more recent pebble games.]
\endroster
\bn
\S3 $\quad$ The canonical example 
\roster
\item "{{}}"  [We deal with random enough graphs and conclude that they
satisfy the 0-1 law so proving the logic cannot express two strong
properties.]
\endroster
\bn
\S4 $\quad$ Closing comments 
\roster
\item "{{}}"  [We present a variant of the criterion (the existence of a
simple $k$-system).  We then define a logic which naturally expresses it.
We comment on defining $N_t[M]$ for ordinals.]
\endroster
\newpage

\head {\S1 The polynomial time logic presented} \endhead  \resetall 
\bigskip

We present below the choiceless polynomial time logic, introduced under the
name $\overset\sim {}\to C PT\text{ime}$ in Blass Gurevich 
Shelah \cite{BGSh:533}.  Knowledge of \cite{BGSh:533} which is phrased with
ASM (abstract state machine) is not required except when we explain how the
definitions fit in \scite{1.3}(4).  See on more relevant works there.
The aim of this logic is to capture statements 
on a (finite) model $M$ in
polynomial time and space \ub{without} arbitrary choices but with no
additional bound on the depth, so its being this logic is a thesis.  
So we are not allowed to use a linear order
on $M$, but if $P^M$ has $\log_2(\|M\|)$ elements we are allowed to list all
subsets of $P^M$, and if e.g. $(|P^M|)! \le \|M\|$ we can list the 
permutations of $P^M$.  Formally for a given $M$, we consider the elements of
$M$ as urelements, and build inductively $N_t = N_t[M]$, with $N_0 = M,
N_{n+1} \subseteq N_t[M] \cup {\Cal P}(N_t[M])$ but the definition is uniform
and $N_t[M]$ should not be too large (i.e. has a polynomial bound) and the
process stops.

Informally, we start with a model $M$ with each element an atom=ure-element,
we successibly define $N_t[M],t$ running on the stages of the ``computations";
to $N_{t+1}[M]$ we add few families of subsets of $N_t$, each of those
defined by a $\psi(-,\bar a)$-formula for some $\bar a$ from $N_t[M]$, and
we update few relations or functions, by defining them from those of the
previous stage.  Those are coded by $c_{t,\ell}$.  We may then check if a
target condition holds, then finishing.  Note that each stage increases the
size of $N_t[M]$ at most by a (fix) power, but in $\|M\|$ steps we can arrive
to a model of size $2^{\|M\|}$.  So we shall have a timing function
$\bold t$ in $\|M\|$, normal polynomial, so when we have wasted too many
resources (e.g. $\|N_t[M]\|+t$) our time is up whether we got an answer or
not.
\bn
More formally
\definition{\stag{1.1} Definition}  1) For a model $M$, with vocabulary
$\tau = \tau_0,\tau$ finite and $\in$ not in $\tau$, let $\tau^+ = \tau_1 = 
\tau \cup \{ \in \}$ considering the elements of $M$ as atoms = urelements, 
we define $V_t[M]$ by induction of $t:V_0[M] = (M_,\in \restriction M)$ 
with $\in \restriction M$ being empty (as we consider the members of $M$ 
as atoms = ``urelements").  Next $V_{t+1}[M]$ is the
model with universe $V_t[M] \cup \{a:a \subseteq V_t[M]\}$ (so we assume
$a \subseteq V_t[M] \Rightarrow a \notin M$ by ``urelements") with the
predicates and individual constants and function symbols of $\tau$ interpreted
as in $M$ (so function symbols in $\tau$ are interpreted as partial functions)
and $\in^{V_{t+1}[M]}$ is $\in \restriction V_{t+1}[M]$. \nl
2) We say $\Upsilon = (\bar \psi,\bar \varphi,\bar c)$ is an inductive scheme
for the language ${\Cal L}_{\text{f.o.}}(\tau^+)$ (where
${\Cal L}_{\text{f.o.}}$ is first order logic) if: letting 
$m_0 = \ell g(\bar \psi),m_1 = \ell g(\bar \varphi)$ we have
\mr
\item "{$(a)$}"  $\bar{\Cal P} = \langle {\Cal P}_\ell:\ell < m_1 \rangle$
is a sequence of unary \footnote{but they will usually be interpreted as
singletons, i.e. acts as individual constants but for the crucial ``times"
$t=0,1$ not necessarily} predicates not in $\tau^+$
\sn
\item "{$(b)$}"  $\bar \psi = \langle \psi_\ell:\ell < m_0 \rangle,
\psi_\ell = \psi_\ell(x;\bar y_\ell)$ is first order in the vocabulary \nl
$\tau_2 = \tau_2[\Upsilon] = \tau_1 \cup \{{\Cal P}_\ell:\ell <
\ell g(\bar c)\}$
\sn
\item "{$(c)$}"  $\bar \varphi = \langle \varphi_\ell:\ell < m_1 \rangle,
\varphi_\ell = \varphi_\ell(x)$ is first order in the vocabulary $\tau_2$.
\ermn
3) For $M,\tau = \tau_0,\tau_1,\Upsilon = (\bar \psi,\bar \varphi,c)$ and
$\tau_2$ as above, we define by induction on $t$ a submodel $N_t = N_t[M] =
N_{\Upsilon,t}[M]$ of $V_t[M]$ and $\bar c_t = \langle c_{t,\ell}:\ell <
m_1 \rangle$ with $c_{t,\ell} \in N_{t+1},N^+_t[M]$ and $P_{t,\ell},
{\Cal P}_{t,k}$ (for $\ell < m_0,k < m_1$) as follows: more exactly
$N_t[M,\Upsilon],c_{t,0}[M,\Upsilon],\bar c_t[M,\Upsilon],P_{t,\ell}[M,
\Upsilon],{\Cal P}_{t,k}[M,\Upsilon]$. 
\sn
We define $N_t[M]$ and ${\Cal P}_{s,\ell},c_{s,\ell}$ for $s < t$ by 
induction on $t$.
\mn
\ub{Case 1 $t=0$}:  $N_t[M] = V_0[M]$.
\mn
\ub{Case 2 $t+1$}:  $N_{t+1}[M]$ is the submodel of $V_{t+1}[M]$ with set of
elements \nl
$N_t[M] \cup \dbcu_{\ell < m_0} {\Cal P}_{t,\ell}[M]$ \nl
where 
\sn
We define ${\Cal P}_{t,k}[M],N^+_t[M]$ and $c_{t,\ell}$ by:

$$
{\Cal P}_{t,k}[M] = \biggl\{ \{a \in N_t[M]:N^+_t[M] \models \psi_k(a,\bar b)
\}:\bar b \in {}^{(\ell g(\bar y_i))}(N_t[M]) \biggr\}
$$

$$
P_{t,\ell}[M] = \{c_{s,\ell}:s < t \text{ and } s \text{ is the immediate
successor of } s\}
$$

$$
N^+_t[M] = (N_t[M],P_{t,0}[M],\dotsc,P_{t,m-1}[M])
$$

$$
c_{t,\ell} = \{a \in N_t[M]:N^+_t[M] \models \varphi_\ell[a]\}
$$
\mn
(so for $t > 1,\psi_k(x,\bar y)$ is actually $\psi'_k(x,\bar y,\bar c_{t-1})
\in {\Cal L}_{\text{f.o.}}(\tau_1)$). 
\mn
\ub{Case 3}:  $t = \infty$ (i.e. $w$) 
\sn
$N_t = \dbcu_{s \le t} N_s$.
\sn
4) We say $\Upsilon$ is standard if some $\psi_i$ guarantees that
$\{s:s < t\} \subseteq N_t[M]$ (remember that we identify the natural number
$t$ with the set $\{0,1,\dotsc,t-1\}$. \nl
5) Let q.d.$(\varphi)$ be the quantifier depth of the formula $\varphi$.
Let $m_{\text{qd}}(\Upsilon) = \text{ Max}\{\text{q.d.}(\varphi):
\varphi \in \{ \varphi_k:k < m_1\} \cup \{\psi_\ell:\ell < m_0\}$ 
and $m_{\text{fv}}(\Upsilon) = \text{ Max}\{\ell g(\bar y_\ell):\ell
< m_0\},\bar m(\Upsilon) = (m_{\text{qd}}(\Upsilon),m_{\text{fv}}
(\Upsilon))\}$.
\nl
6) We may replace above first order logic by a logic ${\Cal L}$.  We let
${\Cal L}_{\text{f.o.}}$ for first order, ${\Cal L}_{\text{card}}$ like 
first order but we demand that $\Upsilon$ is standard and we allow the 
formulas $|\{x:\theta(x,\bar y)\}| = s$.  
We let ${\Cal L}_{\text{card},\bold T}$ (on $\bold T$ see below) be
like ${\Cal L}_{\text{f.o.}}$ but for each $\bold t \in T$ we allow 
the quantifier $(Q_{\bold t}x)\varphi(x,\bar y)$ with $N \models 
(Q_{\bold t}x)\varphi
(x,\bar a)$ iff $\bold t(\text{ure}(N)) < |\{b:N \models \varphi
(b,\bar a]\}|$ where ure$(N)$ is the set of urelements of $N$.  
\enddefinition
\bn
In the definition below the reader can concentrate on $\iota=2$.
\definition{\stag{1.2} Definition}  Let $\bold T$ be a set of functions
$\bold t:\Bbb N \rightarrow \Bbb N \cup \{\infty\}$ and 
${\Cal L}^*$ be a logic
$({\Cal L}_{\text{f.o.}}$ or ${\Cal L}_{\text{card}}$ usually) and let 
$\tau$ be a vocabulary.  If $\bold t$ is constantly $\infty$ we may write
$\infty$.

We define for $\iota = 1,2,3,4$ the logic ${\Cal L}^{\bold T}
_\iota[{\Cal L}^*](\tau)$, for all of those logics the set of sentences for
a vocabulary $\tau$ is a subset of $\Theta = \Theta_\tau = 
\{\theta_{\Upsilon,\chi,\bold t}:
\Upsilon$ an inductive scheme for ${\Cal L}^*(\tau),\chi \in {\Cal L}^*(\tau) 
\text{ and } \bold t \in \bold T\}$, (equal if not said otherwise),
and for most of them we define the stopping time $t_\iota(M,\Upsilon,\bold t)$
or $t_\iota(M,\Upsilon)$ (if $\bold t$ does not matter).  The satisfaction 
relation for ${\Cal L}^{\bold T}_\iota[{\Cal L}^*](\tau)$ is denoted by 
$\models_\iota$.
Also we write $\theta_{\Upsilon,\chi}$ if $\bold t$ does not matter.  We may
let Dom$(\bold t)$ be the set of relevant structures.
\enddefinition
\bn
\ub{Case 1}: For $\iota=1$. 
\sn
$t_\iota[M,\Upsilon] = \text{ Min}\{t:t \ge 2 \text{ and } c_{t,0}[M,\Upsilon]
= c_{t,2}[M,\Upsilon]\}$
\sn
(If there is no such $t \in \Bbb N$ we let it be $\infty$ (i.e. $\omega$ for
set theorists), we could also use ``undefined").
\sn
$M \models_\iota \theta_{\Upsilon,\chi}$ \ub{iff} $N_t[M,\Upsilon] \models
\chi$ for $t = t_\iota[M,\Upsilon]$.
\bn
\ub{Case 2}:  For $\iota=2$. 
\sn
$t_\iota[M,\Upsilon,\bold t] = \text{ Min}\{t:\|N_{t+1}[M,\Upsilon]\| +
(t+1) > \bold t(\|M\|)$ or $c_{t,0}[M,\Upsilon] = c_{t,1}[M,\Upsilon]\}$ and
\mr
\item "{$(a)$}"  if $t = t_\iota[M,\Upsilon,\bold t] < \bold t(\|M\|)$,
\ub{then} the truth value of $\theta_{\Upsilon,\chi,\bold t}$ in $M$ is true
or false \ub{iff} 
$N_t[M,\Upsilon] \models \chi$ or $N_t[M,\Upsilon] \models \neg \chi$
respectively and we write $M \models_\iota \theta_{\Upsilon,\chi,\bold t}$
or $M \models_\iota \neg \theta_{\Upsilon,\chi,\bold t}$ respectively.
\sn
\item "{$(b)$}"  if $t_i[M,\Upsilon,t] = \bold t(\|M\|)$ we say ``$M
\models_\iota \theta_{\Upsilon,\chi,t}$" is undefined and we say ``the truth
value of $\theta_{\Upsilon,\chi,\bold t}$ in $M$" is undefined.
\endroster
\bn
\ub{Case 3}:  $\iota=3$. 
\sn
We restrict ourselves to standard $\Upsilon$ and let
\sn
$t_\iota[M,\Upsilon,\bold t] = \text{ Min}\{t:\|N_{t+1}[M,\Upsilon]\| >
\bold t(\|M\|)\}$ and define $\models_\iota$ as in Case 2.
\bn
\ub{Case 4}:  $\iota=4$. \nl
We restrict ourselves to standard $\Upsilon$ and

$$
\align
t_\iota[M,\Upsilon,\bold t] = \text{ sup}\{t:&\text{for some } k < m_0
\text{ we have } {\Cal P}_{t+1,k}[M,\Upsilon] \\
  &\text{ has } > \bold t(\|M\|) \text{ members}\}
\endalign
$$
\mn
(so it can be $\infty$; i.e. $\omega$ for set theorists, but we can
guarantee ${\Cal P}_{t,0} = \{0,\dotsc,t-1\}$ so that this never happens) and
define $\models_\iota$ as in Case 2. \nl
We can replace in clause (b), $c_0 = c_1$ (i.e. $(\forall x)(P_0(x) \equiv
P_1(x)) \and (\exists x)P_0(x))$ by a sentence $\chi$.
\bn
\ub{\stag{1.3} Discussion}:  1) The most smooth variant for our purpose is
where ${\Cal L}^* = {\Cal L}_{\text{card},\bold T}$ or ${\Cal L}^* =
L_{\text{card},\bold T}$ and $\iota = 4$.  From considering the motivation
the most natural $\bold T$ is $\{n^m:m < \omega\}$, and $\iota =3$. \nl
2) For $\iota = 1,2,3$ some properties of $M$ can be ``incidentally"
expressed by the logic, as the stopping time gives us some information on
concerning cardinality can be expressed.  This suggests preferring the option
$\models_\iota$ undefined in case (b) rather than false. \nl 
3) If you like set theory, you can let $t$ be any ordinal; but this is a
side issue here; see \S4.
\bn
Implicit in \scite{1.2} (and alternative to \scite{1.2}) is
\definition{\stag{1.4} Definition}  Let $M,\Upsilon$ be as in Definition
\scite{1.2} be given. \nl
1) We say $(N,\bar c)$ is an $M$-candidate or $(M,\Upsilon)$-candidate
\mr
\item "{$(a)$}"  $N$ is a transitive submodel of $\dbcu_t V_t[M]$
\sn
\item "{$(b)$}"  $\bar c = \langle c_\ell:\ell < m_1 \rangle,c_\ell \in N$
(or is undefined, so is really a unary relation on $N$ of cardinality
$\le 1$).
\ermn
2) We say $(N',\bar c')$ is the $(\Upsilon,\bold t)$-successor of
$(N,\bar c)$ if $N',c'_\ell$ are defined as in Definition \scite{1.2}, but

$$
|N'| = |N| \cup \dbcu_{\ell < m_1} {\Cal P}_\ell[N,\Upsilon]
$$
\mn
${\Cal P}_\ell[N,\Upsilon] = A_\ell$ is $\bigl\{ \{a:(N,\bar c) \models
\psi_\ell(a,\bar b)\}:\bar b \in {}^{\ell g(\bar g_\ell)}N \bigr\}$ if this
family has $\le \bold t(M)$ members.  $A_\ell$ is empty otherwise. \nl
3) We define $N_t[M,\Upsilon,\bold t]$ and $\bar c_t[M,\Upsilon,\bold t]$
by induction on $t$ as follows:
\mr
\item "{{}}"  for $t=0$ it is $M$ ($c_{t,\ell}$ undefined)
\sn
\item "{{}}"  for $t+1,(N_{t+1},\bar c_{t+1})$ is the $(\Upsilon,\bold t)$-
successor of $(N_t,\bar c_t)$, see below.
\endroster
\enddefinition
\bigskip

\proclaim{\stag{1.4A} Claim}  1) if $(N,\bar c)$ is an 
$(M,\Upsilon)$-candidate, it has at most one $(\Upsilon,\bold t)$-successor.
\nl
2) The various definition fits, so we can use \scite{1.4}.
\endproclaim
\bn
\ub{\stag{1.5} Discussion}:  How do we translate between the definitions above
and \cite{BGSh:533}?
\mr
\widestnumber\item{$(iii)$}
\item "{$(i)$}"  an infinite structure $I$ there corresponds to a
$\tau$-model here
\sn
\item "{$(ii)$}"  a state $A$ there corresponds to a model of the form
$N = N_t[M,\Upsilon]$ in \scite{1.2} and $(N,\bar c)$ in \scite{1.4}
\sn
\item "{$(iii)$}"  dynamic function there corresponds to ${\Cal P}_\iota$
(that is $c_{t,\ell}$) here
\sn
\item "{$(iv)$}"  an object $x$ is active at $A$ in 5.1 there
it corresponds to $x \in N$
\sn
\item "{$(v)$}"  a program \scite{4.7} there corresponds to an $\Upsilon$
in \scite{1.1}(2) here (mainly the first order formulas used);
\sn
\item "{$(vi)$}"  the counting function in \scite{4.8} there corresponds to
the cardinality quantifier (\scite{1.1}(6)) here
\sn
\item "{$(vii)$}"  the polynomial functions $p,q$ in 5.1 there
corresponds to $\bold t \in T$ here.
\ermn
Note that the $c_{t+1,\ell}$ can in the usual set theory manner be actually
$7$-place function from $N_t$ to $N_t$ or $7$-place relation on $N_t$, or be
the universe of $N_t$.  Understanding this to interpret the successor step
there to here we need that all parts of the program are expressible in
${\Cal L}_{\text{f.o.}}$ (or ${\Cal L}_{\text{card}}$).  For the other 
direction we need to show f.o. operations can be expressed by the programs
of ASM there (see 6.1 there), no problem (and not needed to show our
results solved problems there). 
\newpage

\head {\S2 The general elimination of quantifiers
and proofs it's non-expressing} \endhead  \resetall
\bigskip

In Blass Gurevich Shelah \cite{BGSh:533} we deal with the case of 
equality and permutations here we are using partial isomorphisms.  
It seems a reasonably 
precise way, so we shall later, hopefully, get a kind of inverse.
\bn
\ub{\stag{2.0} Discussion}:  Our aim is to have a family ${\Cal F}$ of partial
automorphisms as in Eherenfeucht Fraisse games (or Karp), of the model $M$ 
we analyze, \ub{not} total automorphism (as in \cite{BGSh:533}) - too 
restrictive.  But it has to be lifted to the $N_t$'s.  But their domains (and
ranges) can contain an element of high rank.  So we should not lose anything 
when we get up on $t$.  The solution is $I \subseteq \{A:A \subseteq M\}$ 
closed downward and ${\Cal F}$ (really $\langle {\Cal F}_i:\ell < m_1 
\rangle$), a family of partial automorphisms of $M$.  So every $x \in N$ 
will have a support $A \in I$ and for $f \in {\Cal F}$, its 
action on $A$ determines its action on $x,(G(f)(x)$ in the section notation).
It is not unreasonable
to demand that there is a minimal one, still it is somewhat restrictive (or
we have to add imaginary elements as in \cite{Sh:a} or \cite{Sh:c}, not a 
very appetizing choice). \nl
But how come we in stage $t+1$ succeed to add ``all sets $X = X_{i,h}$"
definable by $\psi_1(x,\bar b),\bar b \in {}^{\ell g(\bar b)}N_t$?

The parameter $b_1,\dotsc,b_m$ each has a support say $A_1,\dotsc,A_m$ all 
in $I$, so we have enough mappings in the family, the new set has in a sense
support $A = \dbcu^m_{\ell =1} A_\ell$, in the sense that suitable partial
mappings do, if $y$ has support $B$ ($BRy$ in this section notation)
$A \cup B \subseteq \text{ Dom}(f),f \restriction A = \text{ id}_A$ the
mapping $f$ induces in $N$, map $y$ to a member of $B$.
\sn
\ub{But} we are not allowed to increase the possible support and $A$ though
a kind of support is probably too large: $I$ is not closed under union.
But, if we add $X = X_{i,\bar b}$ we have to add all similar 
$X' = X_{i,\bar b'}$.  So
our strategy is to say no to looking for a support $A' \in I$.  So fixing
$A'$ we like that if $f \in {\Cal F},f \restriction A' = \text{ id}_{A'},
A \subseteq \text{ Dom}(f)$, then $f \restriction A$ induces a mapping of
$X_{i,\bar b}$ to some $X_{i,\bar b'}$, which we like to demand that will be
equal thus justifying the statement ``$A'$ supports $X$."

How?  We use our bound on the size of the computation.  
\ub{So we need a dichotomy}: 
either there is such $A' \in I$ \ub{or} the number $X_{i,b'}$
defined by $\psi_i(x,\bar b')$ is too large!! \nl
On this dichotomy hangs the proof.
\sn
However, we do not like to state this as a condition on $N_t$ rather on $M$.
We do not ``know" how $\psi_1(x,\bar b')$ will act \ub{but} for any possible
$A'$ this induces an equivalence relation on the images of $A$ (${\Cal F}$ has
to be large enough).

Actually, we can ignore $f$ and develop set theory of elements demanding
support in ${\Cal F}$.  So we break the proof to definition and claims.

We consider three variants of the logic: usual variant to make preservation
clear, and the case with the cardinality quantifier. \nl
We \ub{could use} one ${\Cal F}$ but we use $\langle {\Cal F}_\ell:\ell \le
m_{\text{qd}}(*) \rangle$.  Actually, for much of the treatment only
${\Cal F}_0$ count.
\bn
\ub{Discussion}:  Each $a \in N$ will have support $A \in I$.  Now should we
in $(N,\bar c)$ add the support of each $\bar c_k$ \ub{or} this will be
included?  No!  The $\bar c_t$ will have support $\emptyset$.
\bigskip

\definition{\stag{2.1C} The Main Definition}  1) We say ${\Cal Y} = 
(M,I,{\Cal F})$ is a $k$-system if \nl
$\bar m^* = (m_{\text{qd}}(*),m_{\text{fv}}(*))$ and:
\mr
\item "{$(A)$}"  $I$ is a family of subsets of $|M|$ (the universe of $M$)
closed under subsets and each singleton belongs to it \nl
[hint: intended as first approximation to the possible supports of the
partial automorphism of $M$, the model of course; the intention is that $M$
is finite] \nl
Let $I[m] = \{\dbcu^m_{\ell =1} A_\ell:A_\ell \in I\}$
\sn
\item "{$(B)$}"  ${\Cal F}$ is a family of partial automorphisms of
$M$ such that $A \subseteq \text{ Dom}(f) \and A \in I \Rightarrow f''(A) \in 
I;{\Cal F}$ closed under inverse (i.e. $f \in {\Cal F}_\ell 
\Rightarrow f^{-1} \in {\Cal F}_\ell$) and composition and restriction \nl
\sn
\item "{$(C)$}"  if $f \in {\Cal F}$ then Dom$(f)$ is the union of 
$\le k$ members of $I$
\sn
\item "{$(D)$}"  if $f \in {\Cal F}$ and $A_1,\dotsc,A_{k-1},A_k \in I$ and
$\ell \in \{1,\dotsc,k-1\} \Rightarrow A_\ell \subseteq \text{ Dom}(f)$,
\ub{then} for some $g \in {\Cal F}$ we have 
$$
f \restriction \dbcu^{k-1}_{\ell=1} A_\ell \subseteq g
$$

$$
A_k \subseteq \text{ Dom}(g)
$$
\sn
\item "{$(E)$}"  if $(\alpha)$ then $(\beta)_1 \vee (\beta)_2$ where
{\roster
\itemitem{ $(\alpha)$ }   $h$ is a function from some 
$[m] = \{1,\dotsc,m\}$, Rang$(h)$ belongs to $I[s],E$ is an 
equivalence relation on 
$H_h = \{h':\text{for some } f \in {\Cal F}_0,\text{Rang}(h) \subseteq
\text{ Dom}(f) \text{ and } h' = f \circ h\}$ satisfying \nl
\sn

$\qquad (*) \quad$ \ub{if} $h_1,h_2,h_3,h_4 \in H_h$ and there is \nl

$\qquad \qquad \,\,f \in {\Cal F}$ such that \nl

$\qquad \qquad$ (Rang $h_1)\, \cup$ (Rang $h_3) \subseteq \text{ Dom}(f)$ and
\nl

$\qquad \qquad \,\, f \circ h_1 = h_3,f \circ h_2 = h_2$ \nl

$\qquad \qquad$ \ub{then} $(h_1 E h_2 \Leftrightarrow h_3 E h_4)$
\mn
\itemitem{ $(\beta)_1$ }  there is $u \subseteq [m]$ such that
\nl
Rang$(h \restriction u) \in I$ and \nl
$(\forall h_1,h_2 \in H_h)(h_1 \restriction u = h_2 \restriction u
\rightarrow h_1 E h_2)$
\sn
\itemitem{ $(\beta)_2$ }  the number of $E$-equivalence classes is
$> \bold t(M)$ \nl
[hint: this is how from ``there is some not too large number" we get
``there is a support in $I$"].
\endroster}
\ermn
Note that \wilog \, $h$ is one-to-one. \nl
2) We say $(M,I,{\Cal F})$ is a $(\bold t,\bar m)$-system if $\bar m = (k,s)$
and it is a $(\bold t,s)$-dichotomical $k$-system; we may also say
$(\bold t,k,s)$-system.
\enddefinition
\bigskip

\definition{\stag{2.2B} Definition}  1) We say 
${\Cal Y} = (M,I,{\Cal F})$ is a super $s$-system if it is an
$\bar m^*$-system and in addition
\mr
\item "{$(E)^+$}"   In clause (E), fixing some $A^* \in I[s]$ as a set of 
parameters, the number of equivalence classes is preserved; that is, if
for $\ell = 1,2$ we have $A_\ell \in I[s]$ and for some $m$ and $h_\ell:[m]
\rightarrow A$ with range $\in I$ and ${\Cal E}_{\Cal Y}(A_\ell,h_\ell)$
(see below) and $f \in {\Cal F}$ such that $A_1 \cup \text{ Rang}(h_1)
\subseteq \text{ Dom}(f)$ we have: $f$ maps $A_1$ into $A_2,f \circ h_1 =
h_2$ and $f$ maps $E_1$ to $E_2$ in the natural way, then the number of
$E_1$-equivalence classes is equal to the number of $E_2$-equivalence
classes.
\ermn
2)  Let ${\Cal E}_{\Cal Y}(A,h)$ be defined for 
$A \in I[s]$ and $h$ a function from some $[m]$ into $M$ 
with Rang$(h)$ belonging to $I,A \subseteq \text{ Rang}(h)$ as follows: \nl
${\Cal E}_{\Cal Y}(A,h)$ 
is the set of equivalence relations $E$ on

$$
H_{A,h} = \{h':\text{for some } f \in {\Cal F},f \restriction A =
\text{ id, Rang}(h) \subseteq \text{ Dom}(f) \text{ and } h' = f \circ h\}.
$$
\mn
such that the parallel $(*)$ of clause (E) of Definition \scite{2.1C}(2) 
holds.
(So $H_h$ is replaced by $H_{h,A},f \restriction A$ is the identity).
\nl
Let ${\Cal E}_{\Cal Y}(A) = \dbcu_h {\Cal E}_{\Cal Y}(A,h)$. \nl
3) ${\Cal Y} = (M,I,{\Cal F})$ is a super $(\bold t,(k,s))$-system \ub{if}
it is a super $k$-system and is $(\bold t,s)$-dichotomical.
\enddefinition
\bn
\ub{\stag{2.2C} Observation}:  If $f \in {\Cal F}$ maps $A_1$ to $A_2$, then
it induces a natural extension $\hat f$ of $f$ mapping also ${\Cal E}_{\Cal Y}
(A_1)$ onto ${\Cal E}_{\Cal Y}(A_2)$. 
\bigskip

\definition{\stag{2.3} Definition}  1) Let 
${\Cal Y} = (M,I,{\Cal F})$ be an $k$-system, $M$ a $\tau$-model, $\Upsilon$
is an inductive scheme for ${\Cal L}_{f_0}(\tau^+)$ and $\bar m^* = 
\bar m(\Upsilon)$.

We say that ${\Cal Z} = (N,\bar c,G,R) = (N^{\Cal Z},\bar c^{\Cal Z},
G^{\Cal Z},R^{\Cal Z})$ is a $\Upsilon$-lifting of ${\Cal Y}$ if \mr
\item "{$(a)$}"   $(N,\bar c)$ is an $(M,\Upsilon)$-candidate so $N$ 
is a transitive submodel of set theory with $M$ as its set of urelements
(see \scite{1.4}(1))
\sn
\item "{$(b)$}"  $G$ is a function with domain ${\Cal F}$
\sn
\item "{$(c)$}"  for $f \in {\Cal F},f \subseteq G(f)$ and $G(f)$ is a
partial automorphism of $N$
\sn
\item "{$(d)$}"  if $f \in {\Cal F},g \in {\Cal F},f \subseteq g$ 
then $G(f) \subseteq G(g)$
\sn
\item "{$(e)$}"  $R$ is a two-place relation written $xRy$ such that \nl
$xRy \Rightarrow x \in I \and y \in N$ \nl
[intention: $x$ is a support of $y$]
\sn
\item "{$(f)$}"  if $ARy$ and $f \in {\Cal F},A \subseteq \text{ Dom}(f)$, 
\ub{then} \nl
$y \in \text{ Dom}(G(f))$ and $f \restriction A = \text{ id } \Rightarrow
G(f)(y) = y$
\sn
\item "{$(g)$}"  $(\forall y \in N)(\exists A \in I) ARy$
\sn
\item "{$(h)$}"  if $A \in I$ and $A \subseteq \text{ Dom}(f),y \in
\text{ Dom}(G(f))$ \ub{then} \nl
$ARy \Leftrightarrow f''(A)R(G(f)(y))$
\sn
\item "{$(i)$}"  for $f \in {\Cal F},G(f^{-1}) = (G(f))^{-1}$
\sn
\item "{$(j)$}"  for $f_1,f_2 \in {\Cal F},f = f_2 \circ f_1$
we have $G(f) \subseteq G(f_2) \circ G(f_1)$
\sn
\item "{$(k)$}"  if $c_\ell \in \text{ Dom}(f)$ and $f \in {\Cal F}$ then
$f(c_\ell) = c_\ell$.
\endroster
\enddefinition
\bigskip

\definition{\stag{2.7} Definition}  For ${\Cal Y} = (M,I,{\Cal F})$, a
$k$-system the $0-\Upsilon$-lifting is $(M,\bar c,G,R)$ where
\mr
\item "{$(a)$}"  $G \text{ is the identity on } {\Cal F}$
\sn
\item "{$(b)$}"  $ARy \Leftrightarrow A \in I \and y \in A$
\sn
\item "{$(c)$}"  each $c_\ell$ undefined.
\endroster
\enddefinition
\bn
\ub{\stag{2.8} Fact}:  The $0-\Upsilon$-lifting (in Definition 
\scite{2.7}) exists and is a lifting.
\bigskip

\definition{\stag{2.9} Definition}  Let ${\Cal Y} = (M,I,{\Cal F})$ an
$k$-system and ${\Cal Z} = (N,\bar c,G,R)$ be a $\Upsilon$-lifting. \nl
1) We say $X$ is good or $({\Cal Y},{\Cal Z})$-good if
\mr
\item "{$(a)$}"  $X$ a subset of $N$
\sn
\item "{$(b)$}"  for some $A \in I$ we have $A$ supports $X$ (for our
${\Cal Y}$ and ${\Cal Z}$) which means: \nl
\ub{if} $f \in {\Cal F},BRy$ (so $B \in I,y \in N$) \nl
and $A \cup B \subseteq \text{ Dom}(f)$, and $f \restriction A =
\text{ id}_A$ \nl
\ub{then} $y \in X \Leftrightarrow G(f)(y) \in X$ \nl
(note: $G(f)(y)$ is well defined by clause (g) of \scite{2.3})
\sn
\item "{$(c)$}"  $X \notin N$.
\ermn
2) Let ${\Cal P} = {\Cal P}_{{\Cal Y},{\Cal Z}}$ be the family of 
good subsets of $N$, let ${\Cal R}$ be the two-place relation defined by:
$A{\Cal R}X$ iff $A$ supports $X$, i.e. (b) of part (1) holds. \nl
3) For $f \in {\Cal F}$ and for good $X$ such that $A{\Cal R}X,A \in I$ when
$A \subseteq \text{ Dom}(f)$ we let
\mr
\item "{$(\alpha)$}"   $(G^+(f))(X) = \{y \in N:\text{for some } g \in 
{\Cal F}$ and $y' \in X$ we have \nl
$f \restriction A \subseteq g \text{ and } 
G(g)(y') = y\}$ as $k \ge 3$, if the result is good
\sn
\item "{$(\beta)$}"  $G^+(f) \restriction N = f$.
\ermn
We can prove this, see \scite{2.10} below. \nl
4) We define $E = E_{{\Cal Y},{\Cal Z}}$ as the following relation: 
$X_1 EX_2$ \ub{iff} $X_1,X_2$ are good subsets of $N^{\Cal Z}$ and for 
some $f \in {\Cal F}$ we have $(G^+(f))(X_1) = X_2$; this is an equivalence
relation as $k \ge 2$. \nl
5) ${\Cal Z}' = (N',\bar c',G',R')$ is a successor of ${\Cal Z}$ if both are
$\bar m^*$-systems and:
\mr
\item "{$(a)$}"  $N \subseteq N' \subseteq N \cup {\Cal P}_{{\Cal Y},
{\Cal Z}}$ 
\sn
\item "{$(b)$}"  $X_1 EX_2,X_1 \in {\Cal P}_{{\Cal Y},{\Cal Z}},
X_2 \in {\Cal P}_{{\Cal Y},{\Cal Z}} \Rightarrow [X_1 \in N'
\leftrightarrow X_2 \in N']$
\sn
\item "{$(c)$}"  $G'(f)$ is defined as $G^+(f)$ from part (3) 
\sn
\item "{$(d)$}"  $R'$ is $R \cup [{\Cal R} \restriction (I \times N')]$.
\ermn
6) ${\Cal Z}' = (N',\bar c',G',R')$ is a full $\bold t$-successor of 
${\Cal Z}$ \ub{if} above \nl
$N' = N \cup \{X \in {\Cal P}_{{\Cal Y},{\Cal Z}}:|X/E_{{\Cal Y},{\Cal Z}}| 
\le \bold t(M)\}$.  If we omit $\bold t$ we mean $\bold t(N) = \infty$.  \nl
7) ${\Cal Z}' = (N',\bar c',B',R')$ is the true $(\Upsilon,\bold t)$-successor
of ${\Cal Z}$ if above:
\mr
\item "{$(a)$}"  $(N,\bar c')$ is the $(\Upsilon,\bold t)$-successor of 
$(N,\bar c)$, (see Definition \scite{1.4}(2)).
\ermn
We now prove that Definition \scite{2.9} is O.K.  The functions defined are
functions with the right domain and rang.  The $E$'s are equivalence 
relations.  This is included in the proof of \scite{2.10}.
\enddefinition
\bigskip

\proclaim{\stag{2.10} Claim}  Assume ${\Cal Y}$ is a 
$(\bold t,\bar m^*)$-system (see Definition \scite{2.1C}(2)) and 
${\Cal Y},{\Cal Z},m^*$ are as in Definition \scite{2.9}. \nl
1) The true successor ${\Cal Z}'$ of ${\Cal Z}$ if exists is a successor of 
${\Cal Z}$. \nl
2) A successor ${\Cal Z}'$ of ${\Cal Z}$ is a $\Upsilon$-lifting of 
$({\Cal Y},\bar m^*)$. \nl
3) A full successor of ${\Cal Z}$ exists.
\endproclaim
\bigskip

\demo{Proof}  1) Trivial. \nl
2) We check the clauses in Definition \scite{2.3}.
\enddemo
\bn
\ub{Clause (a)}:  as $N$ is transitive with $M$ its set of ureelements, and
$X \in N' \backslash N \Rightarrow X \in {\Cal P}_N \Rightarrow X \subseteq
N$ also $N'$ is transitive with $M$ its set of urelements.
\bn
\ub{Clause (b)}:  So we have defined $G'$ above.
\bn
\ub{Clause (c)}:  Let $f \in {\Cal F},G'(f) = f'$ and let
$x,y \in N'$ belongs to the domain of $f'$ and we should prove $N' \models 
y \in x \Leftrightarrow N' \models f'(y) \in f'(x)$ (we shall do more 
toward having clause (g) later).
\mn
If $x \in N$, then $f'(x)$ is necessarily in $N$, hence $N' \models y \in
x \Rightarrow y \in N$ and $N' \models z \in f'(x) \Rightarrow z \in N$, so
as $f' \restriction N = f$ we are done.  So we can assume $x \in N' \backslash
N$, so $x$ is a good subset of $N$, so for some $A_0 \in I,A_0Rx$.  We define

$$
\align
z =: \bigl\{ b \in N:&\text{for some } g \in {\Cal F} \text{ and }
b' \in x \text{ we have} \\
  &f \restriction A \subseteq g \text{ and } G(g)(b') = t \bigr\}.
\endalign
$$
\mn
We need the following
\mr
\item "{$(*)_1$}"  $z$ is a good subset of $N$ with $A_1 = f''(A_0)$ a
support of $z$
\sn
\item "{$(*)_2$}"  $x = \{b' \in N:\text{for some } g \in {\Cal F}
\text{ and } b \in z \text{ we have}$ \nl

$\qquad \qquad \quad f^{-1} \restriction 
(f''(A_0)) \subseteq g \text{ and } G(g)(b) = b'\}$
\sn
\item "{$(*)_3$}"  $z$ does not belong to $N$ 
\sn
\item "{$(*)_4$}"  $z = f'(x)$
\sn
\item "{$(*)_5$}"  if $B$ is another support of $x$, then $z'=z$ when
$z' = \{b \in N:\text{ for some } g \in {\Cal F} \text{ and } a \in x$ we
have $f \restriction A \subseteq g$ and $G(g)(a)=b\}$. 
\endroster
\bigskip

\demo{Proof of $(*)_1$}  Suppose $a,b \in N$ and $g_1 \in {\Cal F}$ and
$A_1 \subseteq \text{ Dom}(g_1)$ and \nl
$g_1 \restriction A_1$ is the identity and $a \in \text{ Dom}[G(g_1)]$ and 
$b = G(g_1)(a)$.  Now we should prove that
$a \in z \Leftrightarrow b \in z$.  It is enough to prove $\Rightarrow$ as
applying it to $g^{-1}_1$ we get the other implication.  As $t = G(g_1)(s)$
necessarily by (i) of \scite{2.3} for some support $B_1$ of $s,B_1 \subseteq
\text{ Dom}(g_1)$.

If $a \in z$ then by the definition of $z$ we can find $g \in {\Cal F}$ and
$a' \in X$ such that id$_A \subseteq g$ and $G(g)(a')=a$.  There
is $B_2 \in I$ such that $B_2Ra'$ and $B_2 \subseteq \text{ Dom}(g)$.
As $A,B_1,B_2 \in I$ by clause (D) of Definition \sciteu{2.1} \wilog \,
$B_1 \subseteq \text{ Dom}(g)$ (as $3 \le m_{\text{qd}}(*)$, see
Definition \scite{1.1}).

Let $g' = g_1 \circ g$, so $A \cup B_1 \subseteq \text{ Dom}(g'),
g' \restriction A = g_1 \restriction A = \text{ id}_A$ hence 
$a' \in \text{ Dom}(G(g'))$ and so
$G(g')(a') = G(g_1)(G(g)(a')) = G(g_1)(a) = b$.  (See Definition \scite{2.3},
clause (j).)

So $g',a'$ witness $b \in z$ as required.
\enddemo 
\bigskip

\demo{Proof of $(*)_2$}

Similar.
\enddemo
\bigskip

\demo{Proof of $(*)_3$}  If $z \in N$ there is $A^* \in I$ such that
$A^*Rz$.

Now there is $f_1 \in {\Cal F},f \restriction A \subseteq f_1$ such that
$A^* \subseteq \text{ Rang}(f_1),f_1 \in {\Cal F}$.  
So $z_1 = G(f^{-1}_1)(z)$ is well defined and as in $(*)_2$ the
proof of $(*)_1$ we can check that $\{b \in N:b \in z_1\} = x$; contradiction
to ``$x \notin N$".
\enddemo
\bigskip

\demo{Proof of $(*)_4$}  

Should be clear.
\enddemo
\bn
\demo{Proof of $(*)_5$}  

Should be clear.
\enddemo
\bn
\ub{Clause (d)}:  Check.
\bn
\ub{Clauses (i),(j),(k)}:  Check.
\bn
\ub{Clause (e)}:  See Definition of $R'$.
\bn
\ub{Clause (f)}:  Included in the poof of clause (d).
\bn
\ub{Clause (g)}:  Check.
\bn
\ub{Clause (h)}:

The new case if: $A \subseteq \text{ Dom}(f),A \in I,X \in \text{ Dom}
(G'(f)),X \in N' \backslash N$.  Then $AR'X \Leftrightarrow 
f''(A)R'(G'(f))(X)$ by clause (j) it is enough to prove the implication
$\Rightarrow$.

There is $A_1 \in I$ such that $A_1R'X$.  
If $\neg AR'X$ we can find $g \in {\Cal F},
g \restriction A = \text{ id}_A$, and $z_0 \in \text{ Dom } G(f),z_1 = 
G(f)(z_0)$, such that $z_0 \in X \equiv z_1 \notin X$.  We can find $B_0 \in
I$ such that $A_0Rz_0$ and $B_1 \subseteq \text{ Dom } g$ and let $A_1 =:
g''(A_0)$.  We can find $f_1,f \restriction A \subseteq f_1,B_0 \subseteq
\text{ Dom}(f_1),f_1 \in {\Cal F}$.  Let $B'_0 = f_1(B_0)$.

Now chase arrows.  \nl
3)  Straight.  \hfill$\square_{\scite{2.10}}$
\bigskip

\proclaim{\stag{2.11} Claim}  Assume
\mr
\item "{$(a)$}"  $\Upsilon$ is an inductive scheme for
${\Cal L}_{\text{f.o.}}(\tau^+),\bar m^* = \bar m^*(\Upsilon)$
\sn
\item "{$(b)$}"  $k,s$ satisfy?
\sn
\item "{$(c)$}"  $M$ is a finite $\tau$-model, $\bold t \in \bold T$
\sn
\item "{$(d)$}"  ${\Cal Y} = (M,I,{\Cal F})$ is a $(\bold t,k,s)$-system.
\ermn
\ub{Then}
\mr
\item "{$(\alpha)$}"  \ub{if} $N_t[M,\Upsilon,\bold t],\bar c_t
[M,\Upsilon,\bold t]$ are well defined, \ub{then} for some \nl
${\Cal Z} = (N,\bar c,,R,G)$ a $\Upsilon$-lifting we have 
$(N,\bar c,{\Cal F})(N_t,\bar c_t)$.
\endroster
\endproclaim

\definition{\stag{2.12} Definition}  1) We say that ${\Cal H}$ is a 
$k$-witness to the equivalence of ${\Cal Y}_1$ and ${\Cal Y}_2$ if
\mr
\item "{$(a)$}"  for $\ell =1,2$ we have ${\Cal Y}_\ell = (M_\ell,I_\ell,
{\Cal F}_\ell)$ is a $k$-system
\sn
\item "{$(b)$}"  ${\Cal G}$ is a family of partial isomorphisms from $M_1$
into $M_2$
\sn
\item "{$(c)$}"  for every $g \in {\Cal H}$, we have Dom$(g) \in I_1$,
Rang$(g) \in I_2$
\sn
\item "{$(d)$}"  if $g \in {\Cal H}$ and $f_1 \in {\Cal F}_1$ then $g \circ
f \in {\Cal H}$
\sn
\item "{$(e)$}"  if $g \in {\Cal H}$ and $f_2 \in {\Cal F}_2$ then $f_2 \circ
g \in {\Cal H}$
\sn
\item "{$(f)$}"  if $g \in {\Cal H}$ and $A \in I_1[k-1]$ and $B \in I_1$,
\ub{then} for some $g_1 \in {\Cal H}$ we have $g \restriction A \subseteq
g_1$ and $B \subseteq \text{ Dom}(g_1)$
\sn
\item "{$(g)$}"  if $g \in {\Cal H}$ and $A \in I_2[k-1]$ and $B \in I_2$,
\ub{then} for some $g_1 \in {\Cal H}$ we have $g^{-1} \restriction A
\subseteq g_1$ and $B \subseteq \text{ Rang}(g_1)$.
\ermn
2) We say that ${\Cal H}$ is a $(k,s)$-witness to the equivalence of
$({\Cal Y}_1,\bold t_1)$ and $(\bold Y_2,\bold t_2)$ if
\mr
\widestnumber\item{$(iii)$}  
\item "{$(i)$}"  ${\Cal Y}_\ell$ is $(\bold t,k,s)$-system, i.e.
$(\bold t_\ell,s)$-dichotomical $k$-system for $\ell =1,2$
\sn
\item "{$(ii)$}"  ${\Cal H}$ is a witness to the $k$-equivalence of
${\Cal Y}_1$ and ${\Cal Y}_2$
\sn
\item "{$(iii)$}"  each $g \in {\Cal G}$ preserved the possibility chosen in
the definition or $(\bold t,s)$-dichotomical.
\ermn
3) We say that ${\Cal H}$ is a super $(k,s)$-witness to the equivalence of
$({\Cal Y}_1,\bold t_1)$ and $(\bold Y_2,\bold t_2)$ if 
\mr
\widestnumber\item{$(iii)$}  
\item "{$(i)$}"  ${\Cal Y}_\ell$ is a super $(\bold t_\ell,k,s)$-system
\sn
\item "{$(ii)$}"  ${\Cal H}$ is a witness to the $k$-equivalence of
${\Cal Y}_1$ and ${\Cal Y}_2$
\sn
\item "{$(iii)$}"  each $g \in {\Cal H}$ preserve the cardinalities involved
in the definition of super.
\endroster
\enddefinition
\bn
\ub{\stag{2.13} Main Conclusion}:  Assume
\mr
\item "{$(a)$}"  ${\Cal Y}_\ell = (M_\ell,I_\ell,{\Cal F}_\ell)$ is a
$(\bold t,k,s)$-system for $\ell=1,2$ and $\tau(M_\ell) = \tau$
\sn
\item "{$(b)$}"  ${\Cal H}$ is a $(k,s)$-witness to the equivalence of
$({\Cal Y}_1,\bold t_1)$ and $({\Cal Y}_2,\bold t_2)$
\sn
\item "{$(c)$}"   $\chi \in {\Cal L}_{\text{f.o.}}(\tau^+)$, i.e. a first
order sentence in the vocabulary $\tau^+ = \tau \cup \{\in\}$, and every
subformula of $\chi$ has at most $k$-free variables
\sn
\item "{$(d)$}"  $\Upsilon$ is an inductive scheme not too complicated
relative to $k,s$.
\ermn
\ub{Then}
\mr
\item "{$(\alpha)$}"  the truth value of $\theta_{\Upsilon,\chi,\bold t_1}$
in $M_1$ and $\theta_{\Upsilon,\chi,\bold t_2}$ in $M_2$ are equal except
possibly when: for some $\ell \in \{1,2\}$ we have the truth value of
$\theta_{\Upsilon,\chi,\bold t_\ell}$ in $M_\ell$ is undefined whereas that
of $\theta_{\Upsilon,\chi,\bold t_{3-\ell}}$ in $M_{3-\ell}$ is well defined
and $\bold t[M_\ell,\Upsilon,\bold t_\ell] < t[M_{3-\ell},\Upsilon,
\bold t_{3-\ell}]$
\sn
\item "{$(\beta)$}"  for any 
$t$ if $N^\ell = N_t[M_\ell,\Upsilon,\bold t_\ell]$ is well defined for
$\ell=1,2$, \ub{then} for every sentence $\chi$ such that every subformula
has at most $k$-free variables, we have $N^1 \models \psi \Leftrightarrow
N^2 \models \psi$.
\endroster
\bigskip

\demo{Proof}  Straight.
\enddemo
\bigskip

\proclaim{\stag{2.14} Claim}  In \scite{2.11}, \scite{2.12} we can allow in
the $\psi_i$ and in the $\varphi$ the cardinality quantifier provided that
${\Cal Y}$'s are super.
\endproclaim
\bigskip

\demo{Proof}  Clearer than \scite{2.11}, \scite{2.12}.
\enddemo
\bn
\ub{\stag{2.15} Discussion}  We consider now some variants. \nl
1) We can define a natural equivalence of two $\bar m^*$-systems.  Again the
case with cardinality quantifiers is clearer.

This makes the proof of applications slightly different. \nl
2) We have to consider the stopping times.  If ${\Cal L}^* =
{\Cal L}_{\text{car},\bold T}$ or ${\Cal L}_{\text{card},\bold T}$ this is
natural, (and they are stronger logics than the earlier variants).  
If we still would like to analyze in particular for the others, we 
should be careful how much information can be gotten by the time. \nl
3) We can omit the $c_\ell$'s if the models are rich enough by a first order
formula.  In $N_{t+1}$ reconstruct the sequence $\langle (N_\ell,\bar c_\ell):
\ell \le s \rangle$ (see \S4).
\newpage

\head {\S3 The canonical example} \endhead \resetall
\bn
We apply \S2 to the canonical example: random enough graph.
\definition{\stag{3.1} Definition}  Let $\tau$ be a fixed vocabulary consisting
of predicates only.  We say $M$ is a $(\bold t,k)$-random $\tau$-model 
\ub{if} every quantifier free 1-type over $A \subseteq M,|A| < k$ 
(not explicitly inconsistent) is realized in $M$ by at least 
$\bold s(\|M\|)$ elements.  If $\bold t = \infty$ we may write $k$-random.
\enddefinition
\bigskip

\definition{\stag{3.1A} Definition}  $\bold T_{\text{pol}}$ is $\{f_q:q \in
\Bbb Q,q > 0\}$ where $f_q:\omega \rightarrow \omega$ is $f_q(n) = n^q$ 
or the least integer $\ge n^q$, more exactly.
\enddefinition
\bigskip

\proclaim{\stag{3.2} Claim}  Assume
\mr
\item "{$(a)$}"   $M_\ell$ is $(\bold s_\ell,k)$-random $\tau$-model for
$\ell=1,2$
\sn
\item "{$(b)$}"  $3s \le k$
\sn
\item "{$(c)$}"  $2^{(2^s) \times (\tau)} < \bold t_\ell(M_\ell) <
\bold s_\ell(M_\ell) - s$
\sn
\item "{$(d)$}"  $\Upsilon$ is an inductive scheme for 
${\Cal L}_{\text{f.o.}}(\tau^+)$, $\chi$ a sentence in
${\Cal L}_{\text{f.o.}}(\tau^+)$ and $m_{fv}(\Upsilon) \le s$ and each
subformula of any formula in $\Upsilon$ or $\chi$ has at most $s$ free
variables.
\ermn
\ub{Then} the truth values of $\theta_{\Upsilon,\chi,\bold t_1}$ in $M_1$
and $\theta_{\Upsilon,\chi,\bold t_2}$ in $M_2$ are equal except the case
in \scite{2.12}.
\endproclaim
\bigskip

\demo{Proof}  Let $\ell= 1,2$.  We let $I_\ell = \{A \subseteq M_\ell:|A| \le
q\}$ and 

$$
\align
{\Cal F}_\ell = \{f:&f \text{ is a partial automorphism of } M_\ell \\
  &\text{and Dom}(f) \text{ has at } \le qk \text{ elements}\}
\endalign
$$
\mr
\item "{$(*)_1$}"  ${\Cal Y}_\ell = (M_\ell,I_\ell,{\Cal F}_\ell)$ is a
$k$-system \nl
[why?  the least obvious clause in Definition \scite{2.1C}(1) is clause (D)
which holds by Definition \scite{3.1} above.]
\sn
\item "{$(*)_2$}"  ${\Cal Y}_\ell = (M_\ell,I_\ell,{\Cal F}_\ell)$ is
$(\bold t_\ell,s)$-dichotomical \nl
[why?  let $m \in \Bbb N$ and let $E$ be an equivalence relation on the set
of $h:[m] \rightarrow M_\ell$ satisfying $(*)$ of Definition \scite{2.1C}(2).
Without loss of generality $h$ is one-to-one, so necessarily $m \le s$.  As
$s \le (qk)/2$ there is a quantifier free formula $\varphi(\bar x,\bar y) \in
{\Cal L}(\tau),\ell g(\bar x) = \ell g(\bar y) = m$ such that $h_1 E h_2$ iff
$M_1 \models \varphi[\langle h_1(i):i \in [m] \rangle,\langle h_2[i]:i \in
[m] \rangle$.]
\ermn
Can there be $\bar a_0,\bar a_1 \in {}^m(M_\ell)$ realizing the same
quantifier free type (say $p(\bar x)$) (over the empty set) which are not
$E$-equivalent?  If so we can find $\bar a_2 \in {}^m(M_\ell)$ realizing the
same quantifier free type $p(\bar x)$ and disjoint to $\bar a_0 \char 94
\bar a_1$ (use ``$M_\ell$ is $(3s)$-random"), so \wilog \, $\bar a_0,\bar a_1$
are disjoint.  Now we ask ``are there disjoint $\bar b_0,\bar b_1 \in
{}^m(M_1)$ realizing $p(\bar x)$ which are $E$-equivalent?  If yes, we easily
get a contradiction to ``$E$ an equivalence relation".  So easily there are
at least $s(M_\ell)-m$ pairwise disjoint sequences realizing $p(\bar x)$.
Moreover, easily by transitivity for some $u \subseteq [m]$ we have for
$\bar a,\bar b \in {}^m(M_\ell)$ realizing $p(\bar x),M_\ell \models \varphi
(\bar a,\bar b)$ iff $\bar a \restriction u = \bar b \restriction u$ hence
the number of equivalence classes is $\ge \bold s_\ell(\bold M_\ell)-m$ so we
get one of the allowable answers in Definition \scite{2.1C}(2).
\sn
So assume that there are no $\bar a_0,\bar a_1 \in {}^m(M_\ell)$ realizing the
same quantifier free types over $\emptyset$, hence the number of
$E$-equivalence classes is at most $|S^m_{qf}(\emptyset,M)| \le
2^{(2^m) \times (\tau)}$, which is immaterial for us. but also each 
equivalence class is preserved by any $f \in {\Cal F}$ which gives the other
allowable answer in Definition \scite{2.1C}(2).
\sn
Let

$$
\align
{\Cal H} = \bigl\{ f:&f \text{ is a partial embedding of } M_1 \text{ into }
M_2 \\
  &\text{with \,Dom}(f) \text{ with } \le qk \text{ members} \bigr\}.
\endalign
$$
\mr
\item "{$(*)_3$}"  ${\Cal H}$ is a $(k,s)$-witness to the equivalence of
$({\Cal Y}_1,\bold t_1)$ and $(\bold Y_2,\bold t_2)$ \nl
[why?  straight.]
\ermn
So we can apply \scite{2.12} and get the desired result.
\hfill$\square_{\scite{3.2}}$
\enddemo
\bn
\ub{\stag{3.2A} Conclusion}:  The logic ${\Cal L}^t_4({\Cal L}_{\text{f.o.},
\bold T})(\tau)$ satisfies the 0-1 law for finite random model with a fix
probability for each predicate.
\bn
\ub{\stag{3.2B} Comment}:  We can use time e.g. 
$\|M\|^{\log(\log\|M\|)}$, then in $(*)$, if $M$ has $\mu$-membesr, then
$M$ has to be $k_\mu$-random for appropriate $k_\mu$.
\bigskip

\proclaim{\stag{3.3} Claim}  Consider the vocabulary $\{P\}$, $P$ unary.  For
every sentence $\psi \in {\Cal L}^t_4({\Cal L}_{\text{car},\bold T})(t)$ 
and say time $\bold t$, for any $n$ large enough,
if, $P_1,P_2 \subseteq [n],|P_\ell|,|[n] \backslash P_\ell| \ge
2 \text{ log}_2(\bold t(M))$ then $([n],P_1) \models^{\bold t} \psi
\Leftrightarrow ([n],P_2) \models^{\bold t} \psi$.
\endproclaim
\bigskip

\demo{Proof}  Here we can use automorphisms of $N_t[M,\bold t]$ as in
\cite{BGSh:533} or just use \scite{3.2}.  \hfill$\square_{\scite{3.3}}$
\enddemo
\bigskip

\definition{\stag{3.4} Definition}  1) We say $M$ is a $\tau$-model with
$k$-elimination of quantifiers if for every subset $A_0,A_1$ of $M,|A_0| =
|A_1| < k$ and an isomorphism $f$ from $M \restriction A_0$ onto $M
\restriction A_1$ and $a_0 \in M$ there is $a_1 \in M$ such that $f = f \cup
\{\langle a_0,a_1 \rangle\}$ is an isomorphism from $M \restriction (A_0 \cup
\{a_0\})$ onto $M \restriction (A_1 \cup \{a_1\})$. \nl
2) We replace ``quantifiers" by ``quantifier and counting" if we add: and 
the two sets $\{a'_0 \in M:a'_0,a_0$ realize the same quantifier free type
over $A_0\}$ and $\{a'_1 \in M \backslash a'_1,a_1$ realize the same
quantifier free type over $A_1\}$ has the same number of elements.
\enddefinition
\bigskip

\proclaim{\stag{3.5} Claim}  1) We can in \scite{3.2} weaken the demand
``$M_\ell$ is $(\bold s_\ell,k)$-random $\tau$-model" to
\mr
\item "{$(a)^-_\ell$}"  $(\alpha) \quad M_\ell$ has $k$-elimination of
quantifier
\sn
\item "{${}$}" $(\beta) \quad$ if $\varphi(\bar x,\bar y)$ is a quantifier
free formula defining an equivalence relation \nl

$\quad \,\,\,$ and $\ell g(\bar x) = \ell g(\bar y) \le s$ then the 
number of classes is $\ge \bold t_\ell(M_\ell)$ \ub{or} \nl

$\quad \,\,\,$ each equivalence class is definable by a quantifier free type.
\endroster
\endproclaim
\bigskip

\remark{\stag{3.5A} Remark}  Parallel claims hold for the logic with the
cardinality quantifier.
\endremark
\bigskip

\proclaim{\stag{3.6} Claim}  Choiceless polynomial time + counting logic does
not capture polynomial time.
\endproclaim
\bigskip

\demo{Proof}  Use \scite{2.13} on the question: $|P^M| \ge \|M\|/2$,
similarly to \scite{3.2} with $\tau = \{P\}$.
\enddemo
\newpage

\head {\S4 Closing comments} \endhead  \resetall
\bn
We may consider
\definition{\stag{4.2} Definition}  1) A context is $(K,{\Cal I})$ such that
\mr
\item "{$(a)$}"  $K$ be a class of models with vocabulary (= the set of
predicates) $\tau$
\sn
\item "{$(b)$}"  ${\Cal I}$ is a function such that
\sn
\item "{$(c)$}"  Dom$({\Cal I}) = K$
\sn
\item "{$(d)$}"  ${\Cal I}(M)$ is a family of subsets of $K$, whose union
is $|M|$, and closed under subsets
\sn
\item "{$(e)$}"  ${\Cal I}$ is preserved by isomorphisms.
\ermn
2) In 1) let \nl
Seq$^*_{\Cal I}(M) = \{ \bar a:\bar a \text{ a sequence of
members of } M \text{ of length } \alpha,\text{ Rang}(\bar a) \in {\Cal I}
(M)\}$. \nl
3) We define a logic ${\Cal L}$.  For $k < \omega$ and $\alpha < \omega$ or
just $\alpha$ an ordinal let us define the 
formulas in ${\Cal L}_{k,\alpha}$ by induction on $\alpha$, each 
formula $\varphi$ has the form 
$\varphi(\bar x_0,\bar x_1,\dotsc,\bar x_{k-1}),k_1 \le k$, where the
$\bar x_\ell$'s are pairwise disjoint (finite) sequences of variables and
every variable appearing freely in $\varphi$ appear in one of those
sequences (so any formula is coupled with such $\langle \bar x_0,\dotsc,
\bar x_{k-1} \rangle$, probably some not actually appearing) (we may restrict
to $\bar x_\ell$ finite).
\bn
\ub{$\alpha=0$}:  quantifier free formula; 
i.e. any Boolean combination of atomic ones (with the right variables, of
course).
\bn
\ub{$\alpha +1$}:  $\alpha$ non-limit 
$\varphi(\bar x_0,\dotsc,\bar x_{k_1-1})$ is a Boolean combination of formulas
of the form $(\exists \bar y)\psi(\bar x_{i_0},\dotsc,\bar x_{i_{k_2-2}},
\bar y)$ where $k_2 \le k,\psi(\bar x_{i_0},\dotsc,\bar x_{i_{k_2-2}},\bar y) 
\in {\Cal L}_{k,\alpha}$.
\bn
\ub{$\alpha$ limit}:  ${\Cal L}_{k,\alpha} = \dbcu_{\beta < \alpha}
{\Cal L}_{k,\beta}$.
\bn
\ub{$\alpha +1,\alpha$ limit}:  ${\Cal L}_{k,\alpha +1}$ is the set of
Boolean combinations of members of ${\Cal L}_{k,\alpha}$ of the right
variables.
\mn
Let ${\Cal L}_k = \dbcu_\alpha {\Cal L}_{k,\alpha},{\Cal L}_* =
\dbcu_{k < \omega} {\Cal L}_k,{\Cal L}_{k,< \alpha} = \dbcu_{\beta < \alpha}
{\Cal L}_{k,\beta}$. \nl
4) We now define a satisfaction relation $M \models \varphi(\bar a_0,\dotsc,
\bar a_{k_1-1})$ where $k_1 \le k$ \nl
(depending on ${\Cal I}$). \nl
I.e. we define by induction on $\alpha$, for $\varphi(\bar x_0,\dotsc,
\bar x_{k_1-1}) \in {\Cal L}_{k,\alpha},\bar a_\ell \in 
\text{ Seq}^{\ell g(\bar x_\ell)}_{\Cal I}(M)$, when does $M \models \varphi
[\bar a_0,\dotsc,\bar a_{k_1-1}]$ and when $M \models \neg \varphi[
\bar a_0,\dotsc,\bar a_{k_1-1}]$.  This is done naturally, in particular
$M \models (\exists \bar y)\varphi(\bar a_0,\dotsc,\bar a_{k_2-1},\bar y)$
iff for some $\bar b \in \text{ Seq}^{\ell g(\bar y)}_{\Cal I}(M)$, (so
Rang $\bar b \in {\Cal I}(M)$) we have 
$M \models \varphi[\bar a_0,\dotsc,\bar a_{k_2-1},\bar b]$.
\enddefinition
\bn
\ub{\stag{4.3} Discussion}:  We may replace $M$ by $M^+$, adding 
elements coding each $A \in {\Cal I}(M)$, with decoding by functions, 
still this requires infinitely many functions, we need to actually code 
any sequence listing each $A \in {\Cal I}(M)$.

Still this framework seems to work more smoothly for its purposes.
\bn
\ub{\stag{4.4} Observation}:  In the framework of Definition \scite{4.2},
$M_1 \equiv_{{\Cal L}_k} M_2$ \ub{iff} there is a family ${\Cal F}$ witnessing
it which means
\mr
\item "{$(a)$}"  $f \in {\Cal F} \Rightarrow f$ is a partial isomorphism
from $M_1$ to $M_2$
\sn
\item "{$(b)$}"  $f \in {\Cal F} \Rightarrow \text{ Dom}(f)$ is the union of
$\le k$ members of ${\Cal I}(M_1)$
\sn
\item "{$(c)$}"  $f \in {\Cal F} \Rightarrow \text{ Rang}(f)$ is the union of
$\le k$ members of ${\Cal I}(M_2)$
\sn
\item "{$(d)$}"  if $f \in {\Cal F},A_1,\dotsc,A_{k-1},A_k \in {\Cal I}(M_1),
\ell < k \Rightarrow A_\ell \subseteq \text{ Dom } f$, then for some
$g \in {\Cal F},\dbcu^k_{\ell =1} A_\ell \subseteq \text{ Dom}(g),
g \restriction \dbcu^{k-1}_{\ell =1} A_\ell \subseteq f$
\sn
\item "{$(e)$}"  like (d) for $f^{-1},M_2,M_1$.
\endroster
\bn
\ub{\stag{4.5} Discussion}:  1) In \S2 we can define $N_t[M] \equiv
N_{\Upsilon,t}[M,\bold t]$ for every \ub{ordinal} $t$, and so $V_\Upsilon
[M,\bold t] = \cup\{N_{\Upsilon,\alpha}[M,\bold t]:\alpha 
\text{ an ordinal}\}$, see below.  
Now as in the case $i=4$, the analysis in \S2 works for this but it
is not clear if we can get any interesting things.
\mn
Can this give interesting proofs of consistency for set theory with no choice
but with urelement?
\bigskip

\definition{\stag{4.6} Definition}  1) We say $\Upsilon$ (from Definition
\scite{1.1}) is pure if $m_1[\Upsilon] = 0$ so no $c_\ell$. \nl
2) For pure $\Upsilon$, let ``${\Cal Z}$ is the full $\bold t$-successor of
order $t$ of ${\Cal Z}$" as in \scite{2.9} iterating $t$ times, noting now
the full $\bold t$-successor is unique allowing $t$ to be an ordinal and for
limit ordinal $t$ take just the union.
\enddefinition
\bn
\ub{\stag{4.7} Fact}:  For any $\Upsilon$ we can find $\Upsilon'$ which is
equivalent if we use in Definition \scite{1.2} the case $i=4$ (well when
$\bold t(M_\ell)$ always is $\ge \beth$).  In fact, we can reconstruct the
sequence of $\langle c_{t,\ell}:t' < t \rangle$ in $N_t$.
\bigskip

\definition{\stag{4.8} The Main Definition}  1) We say ${\Cal Y} = 
(M,I,{\Cal F})$ is a $k$-system if \nl
$\bar m^* = (m_{\text{qd}}(*),m_{\text{fv}}(*))$ and:
\mr
\item "{$(A)$}"  $I$ is a family of subsets of $|M|$ (the universe of $M$)
closed under subsets and each singleton belongs to it \nl
[hint: intended as first approximation to the possible supports of the
partial automorphism of $M$, the model of course; the intention is that $M$
is finite] \nl
Let $I[m] = \{\dbcu^m_{\ell =1} A_\ell:A_\ell \in I\}$
\sn
\item "{$(B)$}"  ${\Cal F}$ is a family of partial automorphisms of
$M$ such that $A \subseteq \text{ Dom}(f) \and A \in I \Rightarrow f''(A) \in 
I;{\Cal F}$ closed under inverse (i.e. $f \in {\Cal F}_\ell 
\Rightarrow f^{-1} \in {\Cal F}_\ell$) and composition and restriction \nl
\sn
\item "{$(C)$}"  if $f \in {\Cal F}$ then Dom$(f)$ is the union of 
$\le k$ members of $I$
\sn
\item "{$(D)$}"  if $f \in {\Cal F}$ and $A_1,\dotsc,A_{k-1},A_k \in I$ and
$\ell \in \{1,\dotsc,k-1\} \Rightarrow A_\ell \subseteq \text{ Dom}(f)$,
\ub{then} for some $g \in {\Cal F}$ we have 
$$
f \restriction \dbcu^{k-1}_{\ell=1} A_\ell \subseteq g
$$

$$
A_k \subseteq \text{ Dom}(g)
$$
\sn
\item "{$(E)$}"  if $(\alpha)$ then $(\beta)_1 \vee (\beta)_2$ where
{\roster
\itemitem{ $(\alpha)$ }   $h$ is a function from some 
$[m] = \{1,\dotsc,m\}$, Rang$(h)$ belongs to $I[s],E$ is an 
equivalence relation on 
$H_h = \{h':\text{for some } f \in {\Cal F}_0,\text{Rang}(h) \subseteq
\text{Dom}(f) \text{ and } h' = f \circ h\}$ satisfying \nl
\sn

$\qquad (*) \quad$ \ub{if} $h_1,h_2,h_3,h_4 \in H_h$ and there is \nl

$\qquad \qquad \,\,f \in {\Cal F}$ such that \nl

$\qquad \qquad$ (Rang $h_1)\, \cup$ (Rang $h_3) \subseteq \text{ Dom}(f)$ and
\nl

$\qquad \qquad \,\, f \circ h_1 = h_3,f \circ h_2 = h_2$ \nl

$\qquad \qquad$ \ub{then} $(h_1 E h_2 \Leftrightarrow h_3 E h_4)$
\mn
\itemitem{ $(\beta)_1$ }  there is $u \subseteq [m]$ such that
\nl
Rang$(h \restriction u) \in I$ and \nl
$(\forall h_1,h_2 \in H_h)(h_1 \restriction u = h_2 \restriction u
\rightarrow h_1 E h_2)$
\sn
\itemitem{ $(\beta)_2$ }  the number of $E$-equivalence classes is
$> \bold t(M)$ \nl
[hint: this is how from ``there is some not too large number" we get
``there is a support in $I$"].
\endroster}
\ermn
We can connect this to \S2 as follows.
\enddefinition
\bigskip

\proclaim{\stag{4.10} Claim}  Assume ${\Cal Y}_0 = (M,I,\bar{\Cal F}),
{\Cal L} = (N,\bar c,G),\bar m^*$ are as in \scite{2.1C}, \scite{2.3}.
Then
\mr
\item "{$(a)$}"  if $\varphi(\bar x) \in {\Cal L}^\infty_{m_{\text{qf}}(*)}$
so $\ell g(\bar x) \le m_{\text{qf}}(\bar x)$, and every subformula of
$\varphi(\bar x)$ has $\le m_{\text{qf}}(*)$ free variable, $\bar a \in
{}^{\ell g(\bar x)}N, \dsize \bigwedge_{\ell < \ell g(\bar a)} A_\ell R
a_\ell$ and $\dbcu_{\ell < \ell g(\bar a)} A_\ell \subseteq 
\text{{\rm Dom\/}}(f)$ and $f \in {\Cal F}_0$ then
$$
N \models \varphi[\dotsc,a_\ell,\ldots] \equiv \varphi[\dotsc,G(f)(a_\ell),
\ldots]
$$
\sn
\item "{$(b)$}"  if $\varphi(\bar x)$ has quantifier depth $k$ and
$\le m_{\text{qf}}(*) + m_{\text{qd}}(*) - k$ free variables, 
$f \in {\Cal F}_{m_{\text{qf}}(*)}$ and $\ldots,a_\ell,\dotsc, \in
\text{{\rm Dom\/}}(G(f))$ then 
$$
N \models \varphi[\dotsc,a_\ell,\ldots] \equiv \varphi[\dotsc,G(f)(a_\ell),
\ldots].
$$
\endroster
\endproclaim
\bigskip

\demo{Proof}  Straight.
\enddemo
\bigskip

\ub{\stag{4.11} Conclusion}  1) Assume ${\Cal Y}$ is super (see \scite{2.2B}).
\ub{Then} we can define $R_t,G_t$ for every $t$ ($N_t$ the ``computation"
in time $t$) such that

$$
(M,\bar c_0,G_0,R_0) \text{ is } 0 \text{-lifting}
$$

$$
(N_{t+1},\bar c_{t+1},R_{t+1},\langle G_{t+1} \rangle) \text{ is a lifting,
successor of } (N_t,\bar c_t,G_t,R_t).
$$
\mn
2)  So the formula the $\bar \varphi$ defines is preserved by $f \in
{\Cal F}_0$.
\bigskip

\demo{Proof}  Straight.
\enddemo
\newpage

\shlhetal

\newpage
    
REFERENCES.  
\bibliographystyle{lit-plain}
\bibliography{lista,listb,listx,listf,liste}

\enddocument

\bye